\documentclass{article}
\usepackage{arxiv}
\usepackage[utf8]{inputenc}
\usepackage[T1]{fontenc}
\usepackage{hyperref}
\usepackage{url}
\usepackage{booktabs}
\usepackage{amsfonts}
\usepackage{nicefrac}
\usepackage{microtype}
\usepackage{lipsum}
\usepackage{graphicx}
\usepackage{dcolumn}
\usepackage{bm}
\usepackage{amsmath}
\usepackage{amssymb}
\usepackage{amscd}
\usepackage{amsthm}
\usepackage{amsfonts}
\usepackage{graphicx}
\usepackage{multirow}

\newcommand{\beq}{\begin{equation}}
\newcommand{\eeq}{\end{equation}}
\newcommand{\bea}{\begin{eqnarray}}
\newcommand{\eea}{\end{eqnarray}}
\newcommand{\bc}{\begin{cases}}
\newcommand{\ec}{\end{cases}}

\begin{document}

\date{\today}
\author{ Jos\'{e} M. Amig\'{o} \\
Centro de Investigaci\'{o}n Operativa, \\
Universidad Miguel Hern\'{a}ndez, \\
03202 Elche, Spain \\
\texttt{jm.amigo@umh.es} \And Angel Gim\'{e}nez\\
Centro de Investigaci\'{o}n Operativa, \\
Universidad Miguel Hern\'{a}ndez, \\
03202 Elche, Spain \\
\texttt{a.gimenez@umh.es} }
\title{A simplified algorithm for the topological entropy of multimodal maps}
\maketitle

\begin{abstract}
A numerical algorithm to compute the topological entropy of multimodal maps
is proposed. This algorithm results from a closed formula containing the
so-called min-max symbols, which are closely related to the kneading
symbols. Furthermore, it simplifies a previous algorithm, also based on
min-max symbols, which was originally proposed for twice differentiable
multimodal maps. The new algorithm has been benchmarked against the old one
with a number of multimodal maps, the results being reported in the paper.
In particular, the comparison is favorable to the new algorithm, except in
the unimodal case.
\end{abstract}



\section{Introduction}

\label{sec:1}

Let $f$ be a continuous selfmap of a compact interval $[a,b]\subset \mathbb{R%
}$ with a finite number of turning (or critical) points. Such maps are
generically called multimodal. Then, the topological entropy of $f$ \cite%
{Adler,Walters}, $h(f)$, can be calculated (along other possibilities) with
the formula 
\begin{equation}
h(f)=\lim_{n\rightarrow \infty }\frac{1}{n}\log \ell _{n},  \label{1}
\end{equation}%
where $\ell _{n}$ is shorthand for the \textit{lap number} of $f^{n}$, i.e.,
the number of maximal monotonicity segments of $f^{n}$, the $n$th iterate of 
$f$ \cite{Misiu,Alseda}.

In \cite[Sect. 7]{Amigo2012} a numerical algorithm to compute the
topological entropy of multimodal maps was proposed. Let us point out that
this algorithm generalizes and hence includes a previous one for unimodal
maps published in \cite{Dilao2012}. The algorithm builds on (\ref{1}) by
calculating $\ell_{n}$ with the help of the \textit{min-max symbols} of $f$ 
\cite{Amigo2012,Dilao2012,Dias1,Dilao}, a generalization of the kneading
symbols \cite{Milnor,Melo}. The min-max symbols of a multimodal map not only
locate the iterates of its critical values up to the precision set by the
partition defined by its critical points, as the kneading symbols do, but
they also display their minimum/maximum (or \textquotedblleft
critical\textquotedblright) character. The interesting point is that such an
additional information supposes virtually no extra computational cost.
Indeed, it can be read recursively from a look-up table once the min-max
symbols of the critical values are known.

In this paper we propose a related algorithm which actually approximates the
value of $h(f)$ given by a closed formula involving also the min-max symbols
of $f$. The new algorithm eliminates a formal restriction that, as it turns
out, unnecessarily marred the applicability of the algorithm of \cite%
{Amigo2012}. At the same time, it simplifies the computation scheme of the
latter. We elaborate upon these two points briefly.

With regard to the formal restriction, the theoretical results of \cite%
{Amigo2012} refer to twice differentiable multimodal maps only. However,
numerical simulations with piecewise linear maps of constant slope (and
alternating sign) suggested that the algorithm of \cite{Amigo2012} could be
applied as well to just continuous maps. In this paper we justify the
extension of results from smooth to just continuous maps. Although the proof
turns out to be straightforward, this generalization was not explored in the
previous papers \cite{Amigo2012,Dilao2012} just because these followed the
original approach in \cite{Dias1,Dilao}, which only considered twice
differentiable maps for simplicity.

As for the simplification of the computation scheme, this has to do with the
boundary conditions (or the lack of them). Indeed, the algorithm of \cite%
{Amigo2012} keeps track of the orbits of the boundary points, thus
calculating the exact value of the lap number $\ell_{n}$ in each computation
loop. At the contrary, the new algorithm dispenses with those orbits because
they do not affect the limit (\ref{1}). In fact, Theorem 3 below shows that,
as far as the computation of $h(f)$ is concerned, one may assume that $f$ is
boundary-anchored, i.e., $f(\{a,b\})\subset\{a,b\}$. The result is a compact
expression for the lap number $\ell_{n}$ that makes possible a closed
formula for $h(f)$.

In sum, we fill a theoretical gap in the application of the algorithm in 
\cite{Amigo2012} by showing that continuity of the maps suffices. Moreover,
we abridge the numerical scheme by approximating $\ell_{n}$ in (\ref{1})
with a formula which is exact only for boundary-anchored maps but provides
the right limit \eqref{1} for $h(f)$.

This paper is organized as follows. To make the paper self-contained, we
review in Sect.~2 all the basic concepts, especially the concept of min-max
sequences, needed in the following sections. Most importantly, we extend in
Theorem 1 the transition rules for min-max symbols from twice differentiable
multimodal maps \cite{Amigo2012} to just continuous ones. In Sect.~3 we
introduce some instrumental results which lead in Sect. 4, together with
Theorem 3, to a closed formula for $h(f)$ containing the min-max symbols of $%
f$ (Theorem 4). A formal proof of Theorem 3 has been shifted to the Appendix
in order not to interrupt the flow of ideas. Sect. 5 contains the main
result of the paper, namely, an algorithm for the topological entropy of
(not necessarily smooth) multimodal maps which approximates the value of $%
h(f)$ given in Theorem 4. As way of illustration, this algorithm is put to
test in Sect. 6. First, the new, abridged algorithm is benchmarked in Sects.
6.1 to 6.3 against the full-pledged one, Ref. \cite{Amigo2012}, using smooth
uni-, bi-, and trimodal maps, respectively, borrowed from \cite%
{Amigo2012,Dilao2012}. Finally, in Sect. 6.4, we also compare both
algorithms via piecewise linear, four and five-modal maps of known
topological entropy. It turns out that, except in the unimodal case, the new
algorithm outperforms the old one.

\section{Min-max sequences}

\label{section2}

For the reader's convenience, we use the same notation as in \cite{Amigo2012}
throughout. Let $I$ be a compact interval $[a,b]\subset\mathbb{R}$ and $%
f:I\rightarrow I$ a piecewise monotone continuous map. Such a map is called $%
l$-modal if $f$ has precisely $l$ \textit{turning points} (i.e., points in $%
(a,b)$ where $f$ has a local extremum). Sometimes we speak also of
multimodal maps, in general, or of unimodal maps in the particular case $l=1$%
. Furthermore, assume henceforth that $f$ has local extrema at $%
c_{1}<...<c_{l}$ and is strictly monotone in each of the $l+1$ intervals%
\begin{equation*}
I_{1}=[a,c_{1}),I_{2}=(c_{1},c_{2}),...,I_{l}=(c_{l-1},c_{l}),I_{l+1}=(c_{l},b]%
\text{.}
\end{equation*}
In this case we write $f\in\mathcal{M}_{l}(I).$ When the interval $I$ is
clear from the context or unimportant for the argument, we write just $%
\mathcal{M}_{l}$.

Since the results we obtain below for the calculation of the topological
entropy do not depend on the \textit{shape} of $f$, i.e., on whether $%
f(c_{1} $) is a maximum (positive shape) or a minimum (negative shape), we
assume, unless otherwise stated, that $f$ has positive shape. This implies
that $f(c_{\text{odd}})$ are maxima, whereas $f(c_{\text{even}})$ are
minima, where here and hereafter \textquotedblleft even\textquotedblright\
and \textquotedblleft odd\textquotedblright\ stand for even and odd
subindices, respectively. Hence $f$ is strictly increasing on the intervals $%
I_{\text{odd}}$, and strictly decreasing on the intervals $I_{\text{even}}$.
The points $f(c_{i})$, $1\leq i\leq l$, are called the critical values of $f$
although no differentiability of $f$ at $c_{i}$ is assumed when so doing.

\medskip

\noindent\textbf{Theorem 1}. Let $f\in\mathcal{M}_{l}$ have a positive
shape, and $n\geq1$. Then:

\begin{description}
\item[(a)] If $f^{n}(x)=c_{\text{odd}}$ then $f^{n+1}(x)$ is a maximum. If $%
f^{n}(x)=c_{\text{even}}$ then $f^{n+1}(x)$ is a minimum.

\item[(b)] If $f^{n}(x)$ is a minimum, then 
\begin{equation*}
f^{n+1}(x)\text{ is a }\left\{ 
\begin{array}{cl}
\text{minimum} & \text{if }f^{n}(x)\in I_{\text{odd}}, \\ 
\text{maximum} & \text{if }f^{n}(x)\in I_{\text{even}}.%
\end{array}
\right.
\end{equation*}

\item[(c)] If $f^{n}(x)$ is a maximum, then%
\begin{equation*}
f^{n+1}(x)\text{ is a }\left\{ 
\begin{array}{cl}
\text{maximum} & \text{if }f^{n}(x)\in I_{\text{odd}}, \\ 
\text{minimum} & \text{if }f^{n}(x)\in I_{\text{even}}.%
\end{array}
\right.
\end{equation*}
\end{description}

\medskip

\noindent\textit{Proof}. (a) This is a trivial consequence of $f$ having a
positive shape.

(b) Suppose that $f^{n}(x_{0})$ is a minimum with $f^{n}(x_{0})\in I_{\text{%
odd}}$. Therefore, there exists a neighborhood of $x_{0}$, $U(x_{0})$, such
that $f^{n}(x_{0})\leq f^{n}(x)$ for all $x\in U(x_{0})$. Without
restriction we may assume that $U(x_{0})\subset f^{-n}(I_{\text{odd}})$. It
follows that 
\begin{equation*}
f^{n+1}(x_{0})= f(f^{n}(x_{0}))\leq f(f^{n}(x))= f^{n+1}(x)
\end{equation*}
for all $x\in U(x_{0})$ because $f^{n}(U(x_{0}))\subset I_{\text{odd}}$, an
interval where $f$ is increasing. We conclude that $f^{n+1}(x_{0})$ is a
minimum.

If $f^{n}(x_{0})\in I_{\text{even}}$ then we derive from $f^{n}(x_{0})\leq
f^{n}(x)$ that $f^{n+1}(x_{0})\geq f^{n+1}(x)$ for all $x\in U(x_{0})$
because this time $f^{n}(U(x_{0}))\subset I_{\text{even}}$, an interval
where $f$ is decreasing.

(c) This case follows similarly to (b). $\square$

\medskip

The \textit{itinerary} of $x\in I$ under $f$ is a symbolic sequence 
\begin{equation*}
\mathbf{i}(x)=(i_{0}(x),i_{1}(x),...,i_{n}(x),...)\in%
\{I_{1},c_{1},I_{2},...,c_{l},I_{l+1}\}^{\mathbb{N}_{0}}
\end{equation*}
($\mathbb{N}_{0}\equiv\{0\}\cup\mathbb{N}$), defined as follows:%
\begin{equation*}
i_{n}(x)=\left\{ 
\begin{array}{cl}
I_{j} & \text{if }f^{n}(x)\in I_{j}\text{ }(1\leq j\leq l+1), \\ 
c_{k} & \text{if }f^{n}(x)=c_{k}\text{ }(1\leq k\leq l).%
\end{array}
\right.
\end{equation*}
The itineraries of the critical values,

\begin{equation*}
\gamma^{i}=(\gamma_{1}^{i},...,\gamma_{n}^{i},...)=\mathbf{i}%
(f(c_{i})),\;1\leq i\leq l,
\end{equation*}
are called the \textit{kneading sequences }\cite{Milnor,Melo} of $f$.

\medskip

\noindent\textbf{Definition 1 }\cite{Amigo2012,Dilao2012,Dias1,Dilao}\textbf{%
.} The \textit{min-max sequences} of an $l$-modal map $f$, 
\begin{equation*}
\omega^{i}=(\omega_{1}^{i},\omega_{2}^{i},...,\omega_{n}^{i},...),\;1\leq
i\leq l,
\end{equation*}
are defined as follows:%
\begin{equation*}
\omega_{n}^{i}=\left\{ 
\begin{array}{ll}
m^{\mathbf{\gamma}_{n}^{i}} & \text{if }f^{n}(c_{i})\text{ is a minimum}, \\ 
M^{\mathbf{\gamma}_{n}^{i}} & \text{if }f^{n}(c_{i})\text{ is a maximum.}%
\end{array}
\right.
\end{equation*}
where $\gamma_{n}^{i}$ are kneading symbols.

\medskip

Thus, the \textit{min-max} \textit{symbols} $\omega_{n}^{i}$ have an
exponential-like notation, where the `base' belongs to the alphabet $\{m,M\}$%
, and the `exponent' is a kneading symbol. The extra information of a
min-max symbol $\omega_{n}^{i}$ as compared to a kneading symbol $%
\gamma_{n}^{i}$ is contained, therefore, in the base, which tell us whether $%
f^{n}(c_{i})$ is a minimum ($m$) or a maximum ($M$). Theorem 1 shows that
once the symbol $\omega_{n}^{i}$ of a map with positive shape is known, the
symbol $\omega_{n+1}^{i}$ can be read from the table%
\begin{equation}
\begin{tabular}{|c|c|c|}
\hline
$\omega_{n}^{i}$ & $\rightarrow$ & $\omega_{n+1}^{i}$ \\ \hline\hline
$m^{c_{\text{even}}},M^{c_{\text{even}}}$ & $\rightarrow$ & $m^{\gamma
_{n+1}^{i}}$ \\ \hline
$m^{c_{\text{odd}}},M^{c_{\text{odd}}}$ & $\rightarrow$ & $%
M^{\gamma_{n+1}^{i}}$ \\ \hline
$m^{I_{\text{odd}}},M^{I_{\text{even}}}$ & $\rightarrow$ & $m^{\gamma
_{n+1}^{i}}$ \\ \hline
$m^{I_{\text{even}}},M^{I_{\text{odd}}}$ & $\rightarrow$ & $M^{\gamma
_{n+1}^{i}}$ \\ \hline
\end{tabular}
\ \ \   \label{transition}
\end{equation}
Let us mention for completeness that if $f\in\mathcal{M}_{l}$ has negative
shape, then the transition rules from $\omega_{n}^{i}$ to $\omega_{n+1}^{i}$
read%
\begin{equation}
\begin{tabular}{|c|c|c|}
\hline
$\omega_{n}^{i}$ & $\rightarrow$ & $\omega_{n+1}^{i}$ \\ \hline\hline
$m^{c_{\text{even}}},M^{c_{\text{even}}}$ & $\rightarrow$ & $M^{\gamma
_{n+1}^{i}}$ \\ \hline
$m^{c_{\text{odd}}},M^{c_{\text{odd}}}$ & $\rightarrow$ & $%
m^{\gamma_{n+1}^{i}}$ \\ \hline
$m^{I_{\text{odd}}},M^{I_{\text{even}}}$ & $\rightarrow$ & $M^{\gamma
_{n+1}^{i}}$ \\ \hline
$m^{I_{\text{even}}},M^{I_{\text{odd}}}$ & $\rightarrow$ & $m^{\gamma
_{n+1}^{i}}$ \\ \hline
\end{tabular}
\ \ \ \   \label{transition2}
\end{equation}
instead. This follows \textit{mutatis mutandis} as in the proof of Theorem 1
.

The \textit{transition rules} (\ref{transition}) and (\ref{transition2})
substantiate our claim in the Introduction that, from the point of view of
the computational cost, min-max sequences and kneading sequences are
virtually equivalent.

Therefore, the kneading symbols of $f\in\mathcal{M}_{l}$, along with its 
\textit{initial min-max symbols}, i.e.%
\begin{equation}
\omega_{1}^{i}=\left\{ 
\begin{array}{ll}
M^{\gamma_{1}^{i}} & \text{if }i=1,3,...,2\left\lfloor \frac{l+1}{2}%
\right\rfloor -1, \\ 
m^{\gamma_{1}^{i}} & \text{if }i=2,4,...,2\left\lfloor \frac{l}{2}%
\right\rfloor ,%
\end{array}
\right. \;  \label{initial}
\end{equation}
and the transition rules (\ref{transition}) allow to compute the min-max
sequences of $f\in\mathcal{M}_{l}$ in a recursive way.

A final ingredient (proper of min-max sequences) is the following. Let the $%
i $\textit{th critical line}, $1\leq i\leq l$, be the line $y=c_{i}$ in the
Cartesian product $I\times I$. Min-max symbols split into \textit{bad} and 
\textit{good symbols} with respect to $i$th critical line. Geometrically,
good symbols correspond to local maxima strictly above the line $y=c_{i}$,
or to local minima strictly below the line $y=c_{i}$. All other min-max
symbols are bad by definition with respect to the $i$th critical line. We
use the notation%
\begin{equation*}
\mathcal{B}^{i}=%
\{M^{I_{1}},M^{c_{1}},...,M^{I_{i}},M^{c_{i}},m^{c_{i}},m^{I_{i+1}},...,m^{c_{l}},m^{I_{l+1}}\}
\end{equation*}
for the set of bad symbols of $f\in\mathcal{M}_{l}$ with respect to the $i$%
th critical line. There are $2(l+1)$ bad symbols and $2l$ good symbols with
respect to a given critical line.

Bad symbols appear in all results of \cite{Amigo2012,Dilao2012} concerning
the computation of the topological entropy of $f\in\mathcal{M}_{l}$ via
min-max symbols. In this sense we may say that bad symbols are the hallmark
of this approach.

\section{Auxiliary results}

\label{section3}

Let $s_{\nu}^{i}$, $1\leq i\leq l$, stand for the \textit{number of interior
simple zeros of} $f^{\nu}(x)-c_{i}$, $\nu\geq0$, i.e., solutions of $%
x-c_{i}=0$ ($\nu=0$), or solutions of $f^{\nu}(x)=c_{i}$, $x\in(a,b)$, with $%
f^{\mu}(x)\neq c_{i}$ for $0\leq\mu\leq\nu-1$, and $f^{\nu\prime}(x)\not =0$
($\nu\geq1$). Geometrically $s_{\nu}^{i}$ is the number of \textit{%
transversal} intersections on the Cartesian plane $(x,y)$ of the curve $%
y=f^{\nu}(x)$ and the straight line $y=c_{i}$, over the interval $(a,b)$.
Note that $s_{0}^{i}=1$ for all $i$.

To streamline the notation set%
\begin{equation}
s_{\nu}=\sum_{i=1}^{l}s_{\nu}^{i}  \label{snu}
\end{equation}
for $\nu\geq0$. In particular, 
\begin{equation}
s_{0}=\sum_{i=1}^{l}s_{0}^{i}=\sum_{i=1}^{l}1=l.  \label{seed33}
\end{equation}
According to \cite[Eqn. (31)]{Amigo2012}, the lap number of $f^{n}$, $\ell
_{n}$, satisfies 
\begin{equation}
\ell_{n}=1+\sum\limits_{\nu=0}^{n-1}s_{\nu}=\ell_{n-1}+s_{n-1},
\label{ln-main}
\end{equation}
for $n\geq1$. In particular, $\ell_{1}=\ell_{0}+s_{0}=1+l$.

Furthermore, define%
\begin{equation}
K_{\nu}^{i}=\{(k,\kappa),1\leq k\leq
l,1\leq\kappa\leq\nu:\omega_{\kappa}^{k}\in\mathcal{B}^{i}\},  \label{Kni}
\end{equation}
($\nu\geq1$, $1\leq i\leq l$), that is, $K_{\nu}^{i}$ collects the upper and
lower indices $(k,\kappa)$ of the \textit{bad} symbols with respect to the $%
i $th critical line in all the initial blocks%
\begin{equation*}
\omega_{1}^{1},\omega_{2}^{1},...,\omega_{\nu}^{1};\;\;\omega_{1}^{2},%
\omega_{2}^{2},...,\omega_{\nu}^{2};\;\;...;\;\;\omega_{1}^{l},%
\omega_{2}^{l},...,\omega_{\nu}^{l};
\end{equation*}
of the min-max sequences of $f$. We note for further reference that $K_{\nu
-1}^{i}\subset K_{\nu}^{i}$, the set-theoretical difference being%
\begin{equation}
K_{\nu}^{i}\backslash K_{\nu-1}^{i}=\{(k,\nu),1\leq k\leq
l:\omega_{\nu}^{k}\in\mathcal{B}^{i}\}.  \label{Kni2}
\end{equation}

Finally, set%
\begin{equation}
S_{\nu}^{i}=2\sum_{(k,\kappa)\in K_{\nu}^{i}}s_{\nu-\kappa}^{k},
\label{notation0}
\end{equation}
where $S_{\nu}^{i}=0$ if $K_{\nu}^{i}=\emptyset$, and analogously to (\ref%
{snu}), 
\begin{equation}
S_{\nu}=\sum_{i=1}^{l}S_{\nu}^{i}.  \label{notation}
\end{equation}

We say that $f\in\mathcal{M}_{l}$ is boundary-anchored if $f\{a,b\}\subset
\{a,b\}$. This boundary condition boils down to%
\begin{equation}
f(a)=a,\;\text{and }f(b)=\left\{ 
\begin{array}{ll}
a & \text{if }l\text{ is odd,} \\ 
b & \text{if }l\text{ is even.}%
\end{array}
\right.  \label{boundary}
\end{equation}
for multimodal maps with positive shape, and to%
\begin{equation}
f(a)=b,\;\text{and }f(b)=\left\{ 
\begin{array}{ll}
b & \text{if }l\text{ is odd,} \\ 
a & \text{if }l\text{ is even.}%
\end{array}
\right.  \label{boundary2}
\end{equation}
for multimodal maps with negative shape. As we will see shortly,
boundary-anchored maps have some advantages when calculating the topological
entropy.

\medskip

\noindent\textbf{Theorem 2}. Let $f\in\mathcal{M}_{l}$ be boundary-anchored.
Then%
\begin{equation}
s_{\nu}^{i}=1+\sum_{\mu=0}^{\nu-1}s_{\mu}-S_{\nu}^{i},  \label{sinu}
\end{equation}

\noindent\textit{Proof}. Suppose for the time being that $f$ is twice
differentiable on $[a,b]$ without any restriction at the boundaries. In this
case, it was proved in \cite[Theorem 5.3]{Amigo2012} that%
\begin{equation}
s_{\nu}^{i}=1+\sum_{\mu=0}^{\nu-1}s_{\mu}-S_{\nu}^{i}-\alpha_{\nu}^{i}-%
\beta_{\nu}^{i},  \label{sinu2}
\end{equation}
where $\alpha_{\nu}^{i}$,$\beta_{\nu}^{i}$ are binary variables that vanish
if $f$ is boundary-anchored. Moreover, the relation (\ref{sinu2}) follows
from the transition rules (\ref{transition}) (or (\ref{transition2}) if $f$
has negative shape), which have been proved to hold true also for continuous
multimodal maps in Theorem 1. It follows that (\ref{sinu2}) holds for
continuous, multimodal maps as well. In particular, (\ref{sinu}) holds for
the boundary-anchored ones. $\square$

\medskip

Therefore, the boundary conditions (\ref{boundary}) allow us to express $%
s_{\nu}^{i}$ with the help of some $s_{0}^{k}$, $s_{1}^{k}$, ..., $s_{\nu
-1}^{k}$, $1\leq k\leq l$, via (\ref{sinu}) and (\ref{notation0}). Sum (\ref%
{sinu}) over $i$ from $1$ to $l$ to obtain the compact relation

\begin{equation}
s_{\nu}=l\left( 1+\sum\limits_{\mu=0}^{\nu-1}s_{\mu}\right) -S_{\nu}
\label{account}
\end{equation}
between $s_{0}=l,s_{1},...,s_{\nu}$ and $S_{\nu}$, for all $\nu\geq1$. By (%
\ref{ln-main}) this equation can we rewritten as $s_{\nu}=l\ell_{\nu}-S_{\nu
}$, hence%
\begin{equation}
\ell_{\nu}=\frac{1}{l}(s_{\nu}+S_{\nu}).  \label{account3}
\end{equation}

\medskip

\section{A closed formula for the topological entropy of multimodal maps}

\label{section4}

According to \cite[Lemma 4.4]{Milnor2}, the topological entropy of a
multimodal map depends only on the kneading sequences, i.e., on the
itineraries of the critical values, but not on the itineraries of the
boundary points. This entails that one may assume without restriction the
boundary conditions (\ref{boundary}) or (\ref{boundary2}) when calculating
the topological entropy of $l$-modal maps with positive or negative shape,
respectively. A formal justification is given by the following theorem.

\medskip

\noindent\textbf{Theorem 3.} Let $f\in\mathcal{M}_{l}(I)$. Then there exists 
$F\in\mathcal{M}_{l}(J)$, where $J\supset I$, such that $h(F)=h(f)$ and $F$
is boundary-anchored.

\medskip

See \cite[Lemma 7.7]{Milnor}, and \cite[proof of Lemma 4.4]{Milnor2}. For
the reader's convenience, a proof of Theorem 3 is given in the Appendix.

This being the case, (\ref{1}) and (\ref{account3}) yield the following
result.

\medskip

\noindent\textbf{Theorem 4}. Let $f\in\mathcal{M}_{l}$. Then,%
\begin{align}
h(f) & =\lim_{\nu\rightarrow\infty}\frac{1}{\nu}\log\frac{s_{\nu}+S_{\nu}}{l}
\label{efficient} \\
& =\lim_{\nu\rightarrow\infty}\frac{1}{\nu}\log\frac{1}{l}%
\sum_{i=1}^{l}\left( s_{\nu}^{i}+2\sum_{(k,\kappa)\in
K_{\nu}^{i}}s_{\nu-\kappa}^{k}\right) .  \notag
\end{align}

Eqn. (\ref{efficient}) provides a closed expression for $h(f)$ which
includes the min-max symbols of $f$.

\section{\textbf{A simplified algorithm for the topological entropy}}

\label{section5}

An offshoot of the preceding section is that, when it comes to calculate the
topological entropy of a multimodal map, one can resort to the limit (\ref%
{efficient}), whether the map is boundary-anchored or not. Loosely speaking,%
\begin{equation}
h(f)\simeq\frac{1}{\nu}\log\frac{s_{\nu}+S_{\nu}}{l}  \label{estimat}
\end{equation}
for $\nu$ large enough.

As a matter of fact, the numerical algorithm below estimates $h(f)$ by $%
\frac{1}{\nu}\log\frac{s_{\nu}+S_{\nu}}{l}$ to the desired precision. The
core of the algorithm consists of a loop over $\nu$. Each time the algorithm
enters the loop, the values of $s_{\nu-1}$ and $S_{\nu-1}$ are updated to $%
s_{\nu}$ and $S_{\nu}$, and the current estimation of $h(f)$ is compared to
the previous one. Note that the computation of $S_{\nu}^{i}$, $1\leq i\leq l$%
, requires $s_{0}^{i}=1,s_{1}^{i},...,s_{\nu-1}^{i}$, see (\ref{notation0}),
while the computation of $s_{\nu}^{i}$, $1\leq i\leq l$, requires $%
s_{0}^{i},s_{1}^{i},...,s_{\nu-1}^{i}$, and $S_{\nu}^{i}$, see (\ref{sinu}).

We summarize next the algorithm resulting from (\ref{efficient}) in the
following scheme (\textquotedblleft$A\longrightarrow B$\textquotedblright\
stands for \textquotedblleft$B$ is computed by means of $A$%
\textquotedblright).

\begin{description}
\item[(A1)] \textbf{Parameters: }$l\geq1$ (number of critical points), $%
\varepsilon>0$ (dynamic halt criterion), and $n_{\max}\geq2$ (maximum number
of loops).

\item[(A2)] \textbf{Initialization:} $s_{0}^{i}=1$, and $K_{1}^{i}=\{k,1\leq
k\leq l:\omega_{1}^{k}\in\mathcal{B}^{i}\}$ ($1\leq i\leq l$).

\item[(A3)] \textbf{First iteration: }For $1\leq i\leq l$,\textbf{\ }%
\begin{equation*}
\begin{array}{rcll}
s_{0}^{i},K_{1}^{i} & \longrightarrow & S_{1}^{i},S_{1} & \text{(use (\ref%
{notation0}), (\ref{notation}))} \\ 
s_{0}^{i},S_{1}^{i} & \longrightarrow & s_{1}^{i},s_{1} & \text{(use (\ref%
{sinu}), (\ref{account}))}%
\end{array}%
\end{equation*}

\item[(A4)] \textbf{Computation loop}. For $1\leq i\leq l$ and $\nu\geq2$
keep calculating $K_{\nu}^{i}$, $S_{\nu}^{i}$, and $s_{\nu}^{i}$ according
to the recursions 
\begin{equation}
\begin{array}{rcll}
K_{\nu-1}^{i} & \longrightarrow & K_{\nu}^{i} & \text{(use (\ref{Kni2}), (%
\ref{transition}))} \\ 
s_{0}^{i},s_{1}^{i},...,s_{\nu-1}^{i},K_{\nu}^{i} & \longrightarrow & S_{\nu
}^{i},S_{\nu} & \text{(use (\ref{notation0}), (\ref{notation}))} \\ 
s_{0}^{i},s_{1}^{i},...,s_{\nu-1}^{i},S_{\nu}^{i} & \longrightarrow & s_{\nu
}^{i},s_{\nu} & \text{(use (\ref{sinu}), (\ref{account}))}%
\end{array}
\label{loop}
\end{equation}
until (i) 
\begin{equation}
\left\vert \frac{1}{\nu}\log\frac{s_{\nu}+S_{\nu}}{l}-\frac{1}{\nu-1}\log 
\frac{s_{\nu}+S_{\nu}}{l}\right\vert \leq\varepsilon,  \label{halt}
\end{equation}
or, else, (ii) $\nu=n_{\max}+1$.

\item[(A5)] \textbf{Output.} In case (i) output%
\begin{equation}
h(f)=\frac{1}{\nu}\log\frac{s_{\nu}+S_{\nu}}{l}.  \label{approx}
\end{equation}
In case (ii) output \textquotedblleft Algorithm failed\textquotedblright.
\end{description}

\medskip

As said above, the algorithm (A1)-(A5) simplifies the original algorithm 
\cite{Amigo2012}, which formally consists of the same five steps above but
is based on the exact value of the lap number $\ell_{\nu}$. This entails
that the new algorithm needs more loops to output $h(f)$ with the same
parameter $\varepsilon$ in the halt criterion (\ref{halt}), although this
does not necessarily mean that the overall execution time will be longer
since now less computations are required. In fact, we will find both
situations in the numerical simulations of Sect. 6.

Furthermore, given a halt criterion $\varepsilon$, the execution time
depends as well on the units (i.e., on the base of the logarithm), whichever
algorithm is used. For instance, if logarithms to base $e$ are used (i.e., $%
h(f)$ in \textit{nats}) and $\nu=n_{\text{nat}}$ is the first time that the
halt criterion,%
\begin{equation*}
\left\vert \ln\frac{s_{\nu}+S_{\nu}}{l}-\ln\frac{s_{\nu-1}+S_{\nu-1}}{l}%
\right\vert \leq\varepsilon,
\end{equation*}
happens to hold in the computation loop, then%
\begin{align*}
\left\vert \log_{2}\frac{s_{n_{\text{nat}}}+S_{n_{\text{nat}}}}{l}-\log _{2}%
\frac{s_{n_{\text{nat}}-1}+S_{n_{\text{nat}}-1}}{l}\right\vert & =\frac{1}{%
\ln2}\left\vert \ln\frac{s_{n_{\text{nat}}}+S_{n_{\text{nat}}}}{l}-\ln\frac{%
s_{n_{\text{nat}}-1}+S_{n_{\text{nat}}-1}}{l}\right\vert \\
& \leq\frac{\varepsilon}{\ln2}=1.4427\varepsilon.
\end{align*}
Therefore, if the the halt criterion 
\begin{equation*}
\left\vert \log_{2}\frac{s_{\nu}+S_{\nu}}{l}-\log_{2}\frac{s_{\nu-1}+S_{\nu
-1}}{l}\right\vert \leq\varepsilon
\end{equation*}
for the computation of $h(f)$ in bits does not hold when $\nu=n_{\text{nat}}$%
, i.e.,%
\begin{equation*}
\left\vert \ln\frac{s_{n_{\text{nat}}}+S_{n_{\text{nat}}}}{l}-\ln \frac{%
s_{n_{\text{nat}}-1}+S_{n_{\text{nat}}-1}}{l}\right\vert >(\ln
2)\varepsilon=0.6932\varepsilon\text{,}
\end{equation*}
then the algorithm will not exit the computation loop. We conclude that $n_{%
\text{bit}}\geq n_{\text{nat}}$ with both algorithms, where $n_{\text{bit}}$
is the exit loop when logarithms to base 2 are employed.

Two final remarks:

\begin{description}
\item[R1.] The parameter $\varepsilon$ does not bound the error $\left\vert
h(f)-\frac{1}{\nu}\log\frac{s_{\nu}+S_{\nu}}{l}\right\vert $ but the
difference between two consecutive estimations, see (\ref{halt}). The number
of exact decimal positions of $h(f)$ can be found out by taking different $%
\varepsilon$'s , as we will see in the next section. Equivalently, one can
control how successive decimal positions of $\frac{1}{\nu}\log\frac{s_{\nu
}+S_{\nu}}{l}$ stabilize with growing $\nu$. Moreover, the smaller $h(f)$,
the smaller $\varepsilon$ has to be chosen to achieve a given approximation
precision.

\item[R2.] According to \cite[Thm. 4.2.4]{Alseda}, $\frac{1}{\nu}\log\ell
_{\nu}\geq h(f)$ for any $\nu$. We may expect therefore that the numerical
approximations (\ref{approx}) converge from above to the true value of the
topological entropy with ever more iterations, in spite of the relation $%
\ell_{\nu}=\frac{1}{l}(s_{\nu}+S_{\nu})$ holding in general for
boundary-anchored maps only.
\end{description}

\section{Numerical simulations}

In this section we compute the topological entropy of a variety of
multimodal maps. To this end, a code for arbitrary $l$ was written with
PYTHON, and run on an Intel(R) Core(TM)2 Duo CPU. All the numerical results
will be given with six decimal positions for brevity.

Thus, in Sect. 6.1 to 6.3 we calculate the entropy of families of uni-, bi-,
and trimodal maps, respectively, taken from \cite{Dilao2012} (unimodal case)
and \cite{Amigo2012} (general case). Except for particular values of the
parameters, these maps are not boundary-anchored. The purpose of our choice
is to compare our entropy plots with the plots published in those
references. To complete the picture, we will consider non-smooth maps in
Sect. 6.4. The natural choice are piecewise linear maps of constant slope
because, in this case, the exact value of the topological entropy is known.
In all sections, we are going to compare numerically the performance of the
algorithm presented in Sect. 5 with the general algorithm presented in \cite[%
Sect. 7]{Amigo2012} by means of single maps. For brevity we shall refer to
them as the new algorithm and the old one, respectively.

As for the units, the \textit{nat} is the usual choice in Applied
Mathematics and Physics, while the \textit{bit} is the standard unit in
Information Theory and Communication Technologies. In the following
subsections we are actually going to use both of them despite the fact that,
as shown in Sect. 5, computations with Napierian logarithms are faster to a
given precision. To be specific, we use \textit{bits} in Sect. 6.2, and 6.3
for the sake of comparison with the results published in \cite{Amigo2012},
which are given in that unit.

\subsection{Simulation with 1-modal maps}

Let $\alpha>0$, $-1<\beta\leq0$, and $f_{\alpha,\beta}:[-(1+\beta
),(1+\beta)]\rightarrow\lbrack-(1+\beta),(1+\beta)]$ be defined as \cite[%
Eqn. (29)]{Dilao2012} 
\begin{equation*}
f_{\alpha,\beta}(x)=e^{-\alpha^{2}x^{2}}+\beta.
\end{equation*}
These maps have the peculiarity of showing direct and reverse
period-doubling bifurcations when the parameters are monotonically changed 
\cite[Fig. 3(a)]{Dilao2012}.

Fig. 1 shows the plot of $h(f_{2.8,\beta})$ vs $\beta$ calculated with the
algorithm of Sect. 7. Here $\varepsilon=10^{-4}$ and the parameter $\beta$
was increased in steps of $\Delta\beta=0.001$ from $\beta=-0.999$ to $%
\beta=0 $. Upon comparing Fig. 1 with Fig. 3(b) of \cite{Dilao2012}, we see
that both plots coincide visually, except for the two vanishing entropy
tails. We conclude that $\varepsilon=10^{-4}$ is not small enough to obtain
reliable estimations of the topological entropy for vanishing values of $%
h(f_{2.8,\beta })$. This fact can also be ascertained numerically by taking
different values of $\varepsilon$, as we do in the table below.

\begin{figure}[h]
\centering
\includegraphics[width=10cm]{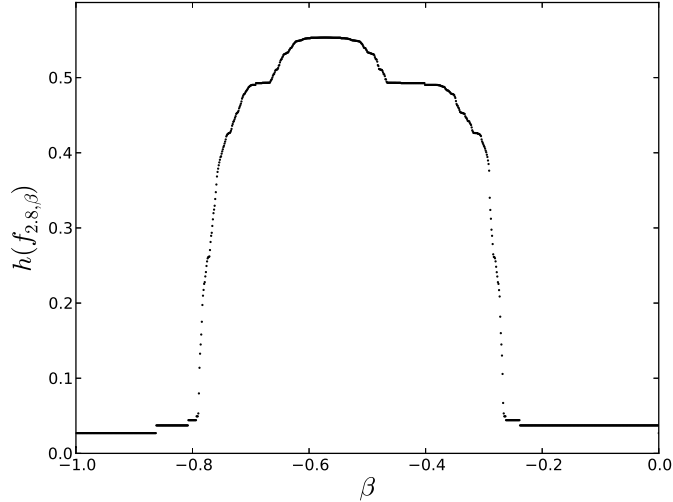}
\caption{Plot of $h(f_{2.8,\protect\beta})$ in nats vs $\protect\beta$, $-1<%
\protect\beta\leq0$ ($\protect\varepsilon=10^{-4},\Delta\protect\beta=0.001$%
).}
\label{fig1a}
\end{figure}

To compare the convergence speed and execution time of the old and the new
algorithm, we have computed $h(f_{2.8,-0.5})$ with both algorithms for
different $\varepsilon$'s. The number of loops $n$ needed to achieve the
halt condition $\varepsilon=10^{-d}$, $4\leq d\leq7$, and the execution time 
$t$ (in seconds) are listed in Table 1. The columns $h_{old}$, $n_{old}$,
and $t_{old}$ were obtained with the old algorithm, while the columns $%
h_{new}$, $n_{new}$, and $t_{new}$ were obtained with the new one. For $%
\varepsilon =10^{-4}$ it exceptionally holds $t_{old}>t_{new}$, otherwise $%
t_{old}<t_{new}$. Furthermore, we conclude from Table 1 that $h(f)=0.52...$
nats with either algorithm and $\varepsilon=10^{-6}$, both decimal digits
being exact. If $\varepsilon=10^{-7}$ the old algorithm fixes the third
decimal digit, $h(f)=0.524...$ nats, whereas the new algorithm does not.

\begin{table}[]
\centering
\begin{tabular}{|c|c|c|c|c|c|c|}
\hline
& $h_{old}$ & $n_{old}$ & $t_{old}$ & $h_{new}$ & $n_{new}$ & $t_{new} $ \\ 
\hline
$\varepsilon=10^{-4}$ & 0.531968 & 81 & 0.031076 & 0.534106 & 101 & 0.021248
\\ \hline
$\varepsilon=10^{-5}$ & 0.526645 & 253 & 0.179558 & 0.527305 & 318 & 0.193149
\\ \hline
$\varepsilon=10^{-6}$ & 0.524935 & 797 & 1.684213 & 0.525142 & 1004 & 
1.912784 \\ \hline
$\varepsilon=10^{-7}$ & 0.524391 & 2519 & 16.369158 & 0.524456 & 3174 & 
18.900032 \\ \hline
\end{tabular}%
\caption{Comparison of performances when computing $h(f_{2.8,-0.5})$ in
nats. }
\label{table1}
\end{table}

Fig. 2 depicts the values of $h(f_{\alpha,\beta})$ for $2\leq\alpha\leq3$, $%
-1<\beta\leq0$, $\varepsilon=10^{-4}$, and $\Delta\alpha,\Delta\beta=0.01$.

\begin{figure}[h]
\centering
\includegraphics[width=10cm]{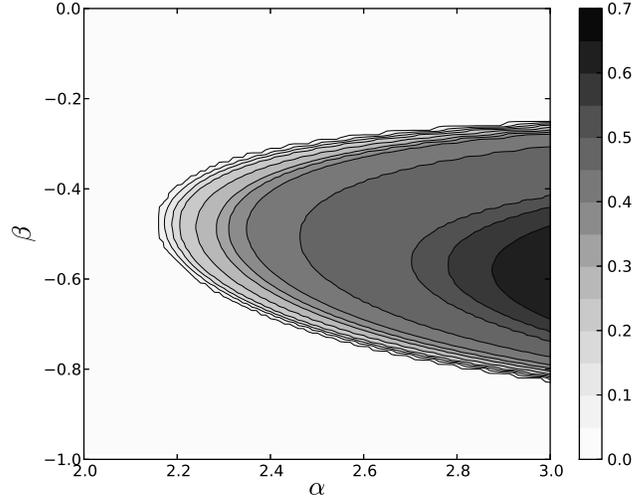}
\caption{Level sets of $h(f_{\protect\alpha ,\protect\beta})$ in nats vs $%
\protect\alpha,\protect\beta$, $2\leq\protect\alpha\leq3$, and $-1<\protect%
\beta\leq0$ ($\protect\varepsilon=10^{-4},\Delta\protect\alpha=\Delta\protect%
\beta=0.01$).}
\label{fig1b}
\end{figure}

\subsection{Simulation with 2-modal maps}

Let $0\leq v_{2}<v_{1}\leq1$ and $f_{v_{1},v_{2}}:[0,1]\rightarrow%
\lbrack0,1] $ be defined as \cite[Sect. 8.1]{Amigo2012} 
\begin{equation*}
f_{v_{1},v_{2}}(x)=(v_{1}-v_{2})(16x^{3}-24x^{2}+9x)+v_{2},
\end{equation*}
These maps have convenient properties for numerical simulations as they
share the same fixed critical points,%
\begin{equation*}
c_{1}=1/4,\;c_{2}=3/4,
\end{equation*}
the critical values are precisely the parameters,%
\begin{equation*}
f_{v_{1},v_{2}}(1/4)=v_{1},\;f_{v_{1},v_{2}}(3/4)=v_{2},
\end{equation*}
and the values of $f$ at the endpoints are explicitly given by the
parameters as follows:%
\begin{equation*}
f_{v_{1},v_{2}}(0)=v_{2},\;f_{v_{1},v_{2}}(1)=v_{1}.
\end{equation*}

Fig. 3 shows the plot of $h(f_{1,v_{2}})$ vs $v_{2}$, $0\leq v_{2}<1$,
computed with the new algorithm, $\varepsilon=10^{-4}$, and $\Delta
v_{2}=0.001$. Again, this plot coincides visually with the same plot
computed with the old algorithm \cite[Fig. 4]{Amigo2012} except for the
vanishing entropy tail, which indicates that $\varepsilon=10^{-4}$ is too
large a value for obtaining accurate estimates in that parametric region. 
\begin{figure}[h]
\centering
\includegraphics[width=10cm]{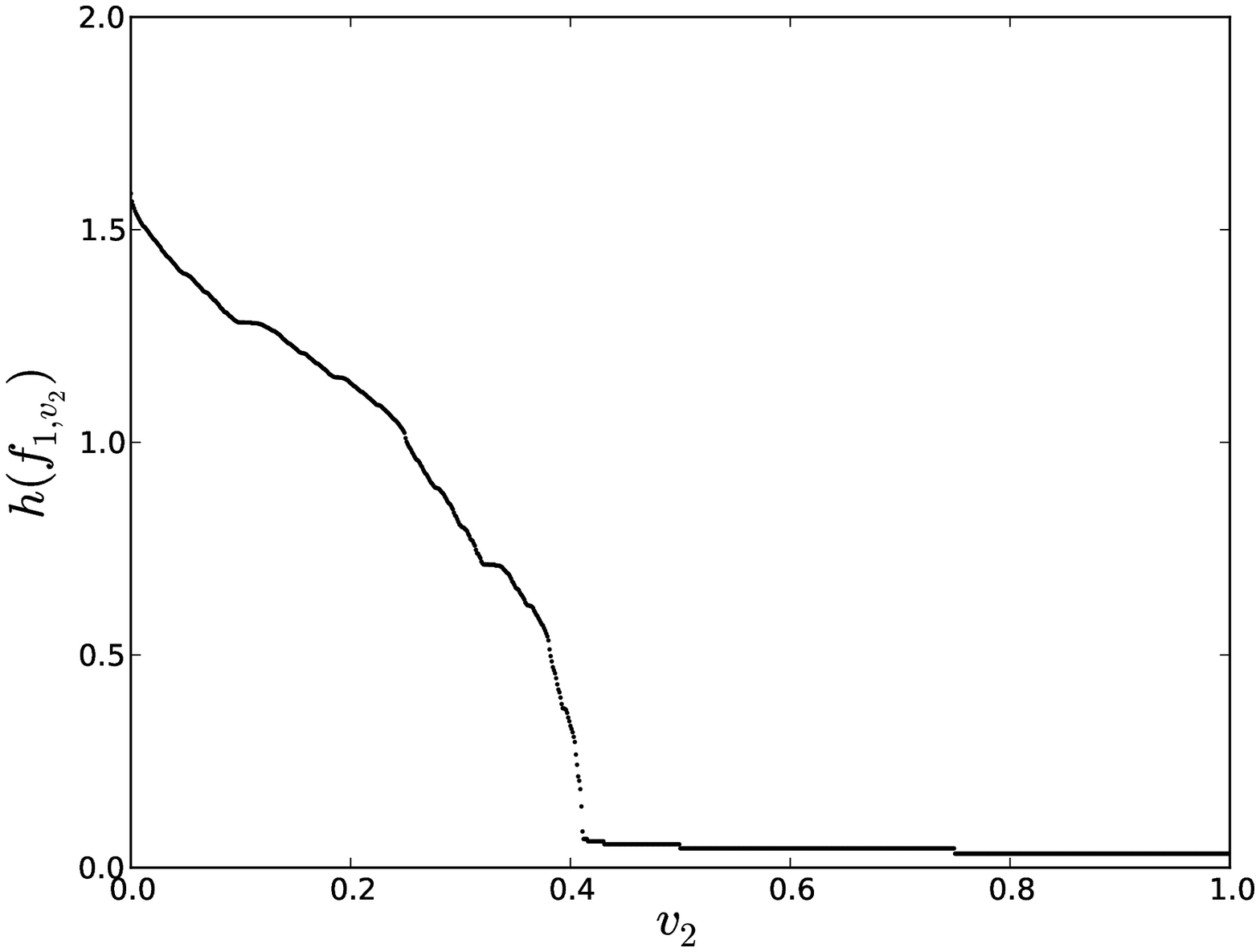}
\caption{Plot of $h(f_{1,v_{2}})$ in bits vs $v_{2}$, $0\leq v_{2}\leq1$ ($%
\protect\varepsilon=10^{-4},\Delta v_{2}=0.001$).}
\end{figure}

Table 2 displays the performance of the new algorithm as compared to the old
one when computing $h(f_{0.9,0.1})$. This time $t_{old}>t_{new}$ for $%
\varepsilon=10^{-d}$, $4\leq d\leq7$ (as in Table 1). Furthermore, we obtain
two correct decimal digits of the topological entropy, $%
h(f_{0.9,0.1})=0.60...$ bits, with both algorithms and $\varepsilon=10^{-6}$.

\begin{table}[h]
\centering
\begin{tabular}{|c|c|c|c|c|c|c|}
\hline
& $h_{old}$ & $n_{old}$ & $t_{old}$ & $h_{new}$ & $n_{new}$ & $t_{new}$ \\ 
\hline
$\varepsilon=10^{-4}$ & 0.619682 & 195 & 0.286922 & 0.622100 & 218 & 0.253133
\\ \hline
$\varepsilon=10^{-5}$ & 0.606568 & 613 & 2.665108 & 0.607310 & 688 & 2.485049
\\ \hline
$\varepsilon=10^{-6}$ & 0.602385 & 1938 & 26.238006 & 0.602622 & 2173 & 
24.890648 \\ \hline
$\varepsilon=10^{-7}$ & 0.601062 & 6125 & 271.074381 & 0.601137 & 6871 & 
265.198039 \\ \hline
\end{tabular}%
\caption{Comparison of performances when computing $h(f_{0.9,0.1})$ in bits.}
\label{table2}
\end{table}

Fig. 4 depicts the values of $h(f_{v_{1},v_{2}})$ for $0\leq v_{2}\leq
v_{1}-0.5$, $\varepsilon=10^{-4}$, and $\Delta v_{1},\Delta v_{2}=0.01$.

\begin{figure}[h]
\centering
\includegraphics[width=10cm]{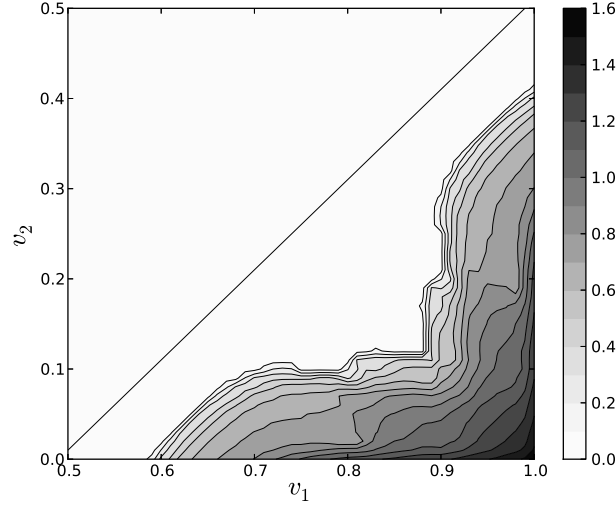}
\caption{Level sets of $h(f_{v_{1},v_{2}})$ in bits vs $v_{1},v_{2}$, $0\leq
v_{2}\leq v_{1}-0.5$ ($\protect\varepsilon =10^{-4},\Delta v_{1}=\Delta
v_{2}=0.01$).}
\label{fig2b}
\end{figure}

\subsection{Simulation with 3-modal maps}

Consider next the $3$-modal maps $f_{v_{2},v_{3}}:[0,1]\rightarrow%
\lbrack0,1] $ defined by the quartic polynomials \cite[Sect. 8.2]{Amigo2012}%
\begin{align*}
f_{v_{2},v_{3}}(x) & =\frac{4\left[ \left( 2\sqrt{2}-1\right) v_{2}-2v_{3}%
\right] x}{2(2\sqrt{2}+1)v_{3}-7v_{2}}\left[ 4\left( 1+2\sqrt{2}\right)
(x-1)(1-2x)^{2}v_{3}\right. \\
& +\left. \left( -56x^{3}+20\left( 4+\sqrt{2}\right) x^{2}-\left( 37+18\sqrt{%
2}\right) x+3\sqrt{2}+5\right) v_{2}\right] ,
\end{align*}
where $0\leq v_{2}<v_{3}\leq1$. The critical points of $f_{v_{2},v_{3}}$ are 
\begin{equation*}
c_{1}=\frac{-\sqrt{2}\text{$v_{2}$}-4\text{$v_{2}$}+12\sqrt{2}\text{$v_{3}$}%
-8\text{$v_{3}$}}{8\left( -7\text{$v_{2}$}+4\sqrt{2}\text{$v_{3}$}+2\text{$%
v_{3}$}\right) },\quad c_{2}=1/2,\quad c_{3}=\frac{1}{4}(2+\sqrt {2}).
\end{equation*}
Moreover this family verifies $f_{v_{2},v_{3}}(0)=0$, $%
f_{v_{2},v_{3}}(c_{2})=v_{2}$, $f(c_{3})=v_{3}$, and

\begin{equation*}
f_{v_{2},v_{3}}(1)=\frac{4\left( 5\sqrt{2}-8\right) \text{$v_{2}$}\left(
\left( 2\sqrt{2}-1\right) \text{$v_{2}$}-2\text{$v_{3}$}\right) }{-7\text{$%
v_{2}$}+4\sqrt{2}\text{$v_{3}$}+2\text{$v_{3}$}}.
\end{equation*}

Fig. 5 shows the plot of $h(f_{v_{2},1})$ vs $v_{2}$, $0\leq v_{2}<1$,
computed with the new algorithm, $\varepsilon=10^{-4}$, and $\Delta
v_{2}=0.001$. Once more, this plot coincides visually with the same plot
computed with the old algorithm \cite[Fig. 7 (left)]{Amigo2012} except for
the vanishing entropy tail, which again indicates that $\varepsilon=10^{-4}$
is too large a value for obtaining accurate estimates in that parametric
region.

\begin{figure}[h]
\centering
\includegraphics[width=10cm]{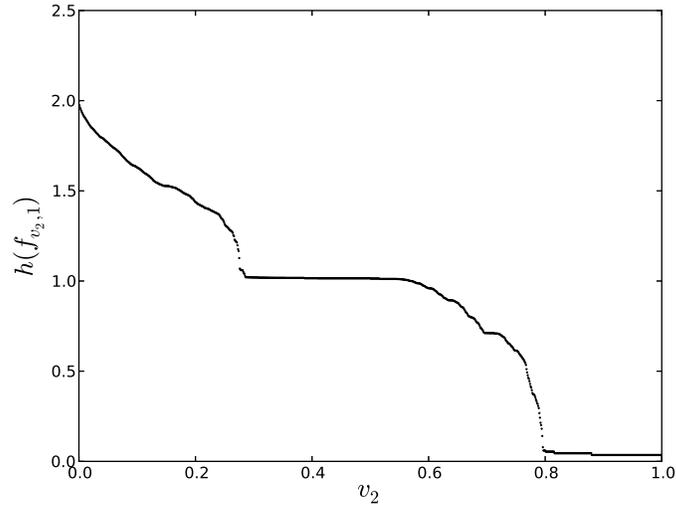}
\caption{Plot of $h(f_{v_{2},1})$ in bits vs $v_{2}$, $0\leq v_{2}<1$ ($%
\protect\varepsilon=10^{-4},\Delta v_{2}=0.001$).}
\label{fig3a}
\end{figure}

Table 3 displays the performance of the new algorithm as compared to the old
one when computing $h(f_{0.7,1})$. Also this time $t_{old}>t_{new}$ for $%
\varepsilon=10^{-d}$, $4\leq d\leq7$ (as in Table 1 and 2). Furthermore, we
obtain two correct decimal digits of the topological entropy, $%
h(f_{0.7,1})=0.69...$ bits, with both algorithms and $\varepsilon=10^{-6}$.

\begin{table}[h]
\centering
\begin{tabular}{|c|c|c|c|c|c|c|}
\hline
& $h_{old}$ & $n_{old}$ & $t_{old}$ & $h_{new}$ & $n_{new}$ & $t_{new}$ \\ 
\hline
$\varepsilon=10^{-4}$ & 0.710321 & 162 & 0.422875 & 0.711709 & 177 & 0.351586
\\ \hline
$\varepsilon=10^{-5}$ & 0.699339 & 511 & 4.028703 & 0.699793 & 557 & 3.448455
\\ \hline
$\varepsilon=10^{-6}$ & 0.695855 & 1615 & 39.882275 & 0.696000 & 1759 & 
34.853025 \\ \hline
$\varepsilon=10^{-7}$ & 0.694752 & 5105 & 417.489291 & 0.694798 & 5561 & 
368.900286 \\ \hline
\end{tabular}%
\caption{Comparison of performances when computing $h(f_{0.7,1})$ in bits.}
\label{table3}
\end{table}

Fig. 6 depicts the values of $h(f_{v_{2},v_{3}})$ for $v_{2}+0.3\leq
v_{3}\leq1$, $\varepsilon=10^{-4}$, and $\Delta v_{2},\Delta v_{3}=0.01$.

\begin{figure}[h]
\centering
\includegraphics[width=10cm]{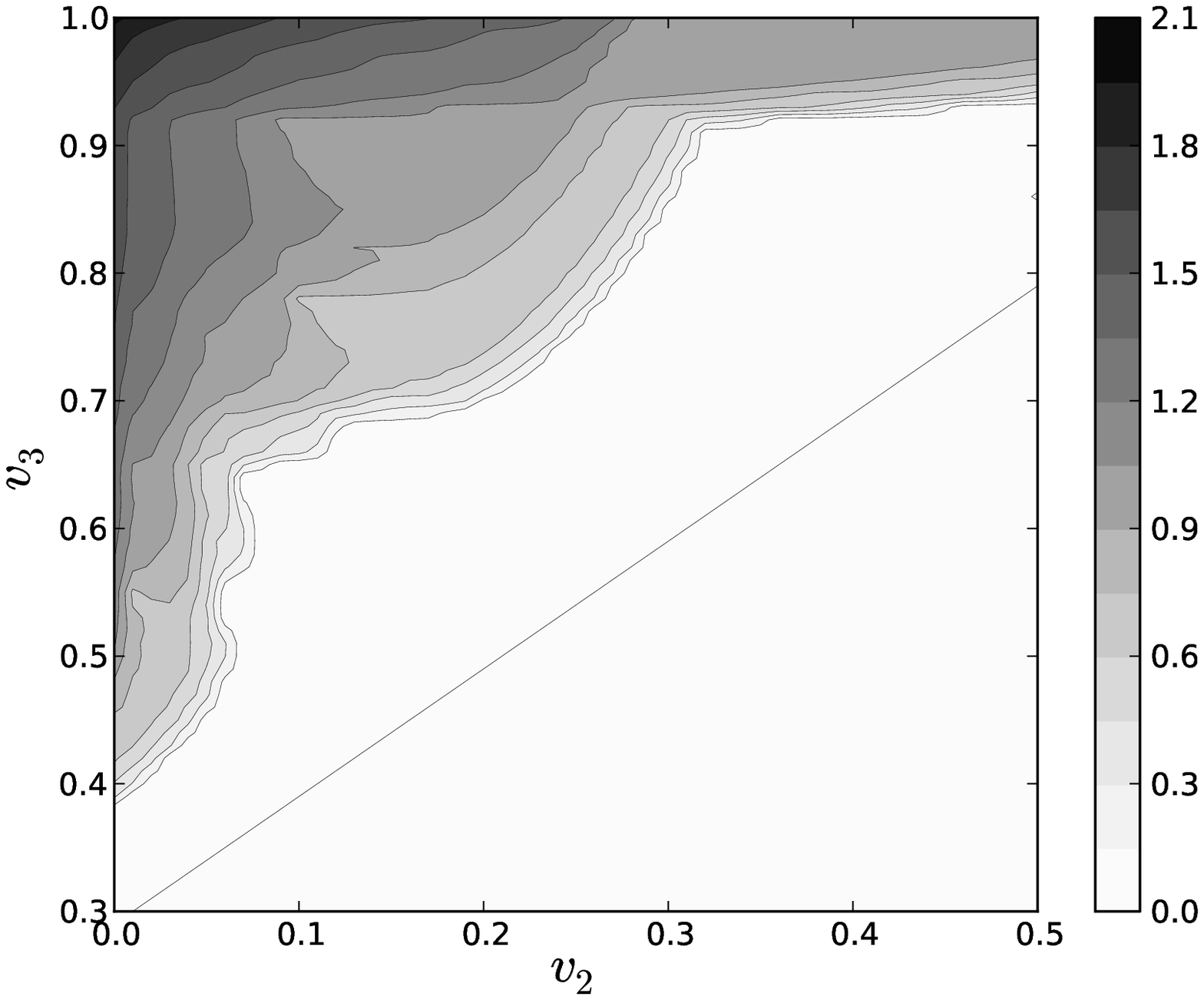}
\caption{Level sets of $h(f_{v_{2},v_{3}})$ in bits vs $v_{2},v_{3}$, $%
v_{2}+0.3\leq v_{3}\leq1$ ($\protect\varepsilon =10^{-4},\Delta v_{1}=\Delta
v_{2}=0.01$).}
\label{fig3b}
\end{figure}

\subsection{Simulation with higher multimodal maps}

Last but not least, we are going to compare the performance of the old and
new algorithms with the 4- and 5-modal maps of Fig. 7. These are piecewise
linear maps on $[0,1]$, with constant slope $s=\pm1.5$, critical points%
\begin{equation*}
c_{1}=\tfrac{3}{10},\;c_{2}=\tfrac{23}{60},\;c_{3}=\tfrac{7}{15},\;c_{4}=%
\tfrac{11}{20},
\end{equation*}
and critical values%
\begin{equation*}
f(c_{1})=f(c_{3})=0.450,\;f(c_{2})=f(c_{4})=0.325,
\end{equation*}
in the $l=4$ case, while%
\begin{equation*}
c_{1}=0.3,\;c_{2}=0.4,\;c_{3}=0.5,\;c_{4}=0.6,\;c_{5}=0.7,
\end{equation*}
and%
\begin{equation*}
f(c_{1})=f(c_{3})=f(c_{5})=0.45,\;f(c_{2})=f(c_{4})=0.30,
\end{equation*}
in the $l=5$ case. By \cite[Corollary 4.3.13]{Alseda}, 
\begin{equation*}
h(f)=\max\{0,\ln\left\vert s\right\vert \}=\ln1.5=0.40547\text{ nats}
\end{equation*}
in both cases. At variance with the previous examples in Sects. 6.1 to 6.3,
these two maps are non-smooth and boundary anchored.

\begin{figure}[h]
\centering
\includegraphics[width=6cm]{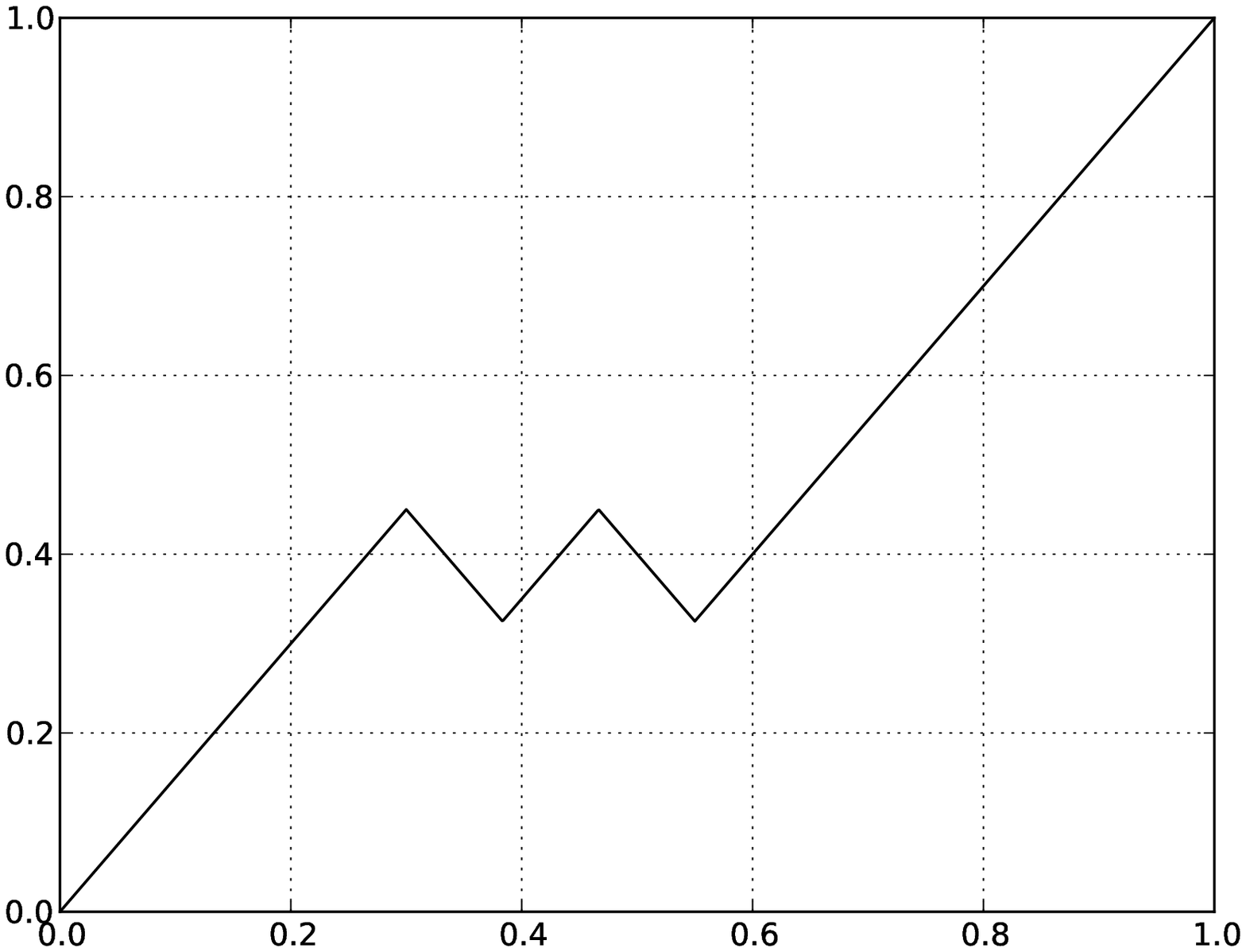} \includegraphics[width=6cm]{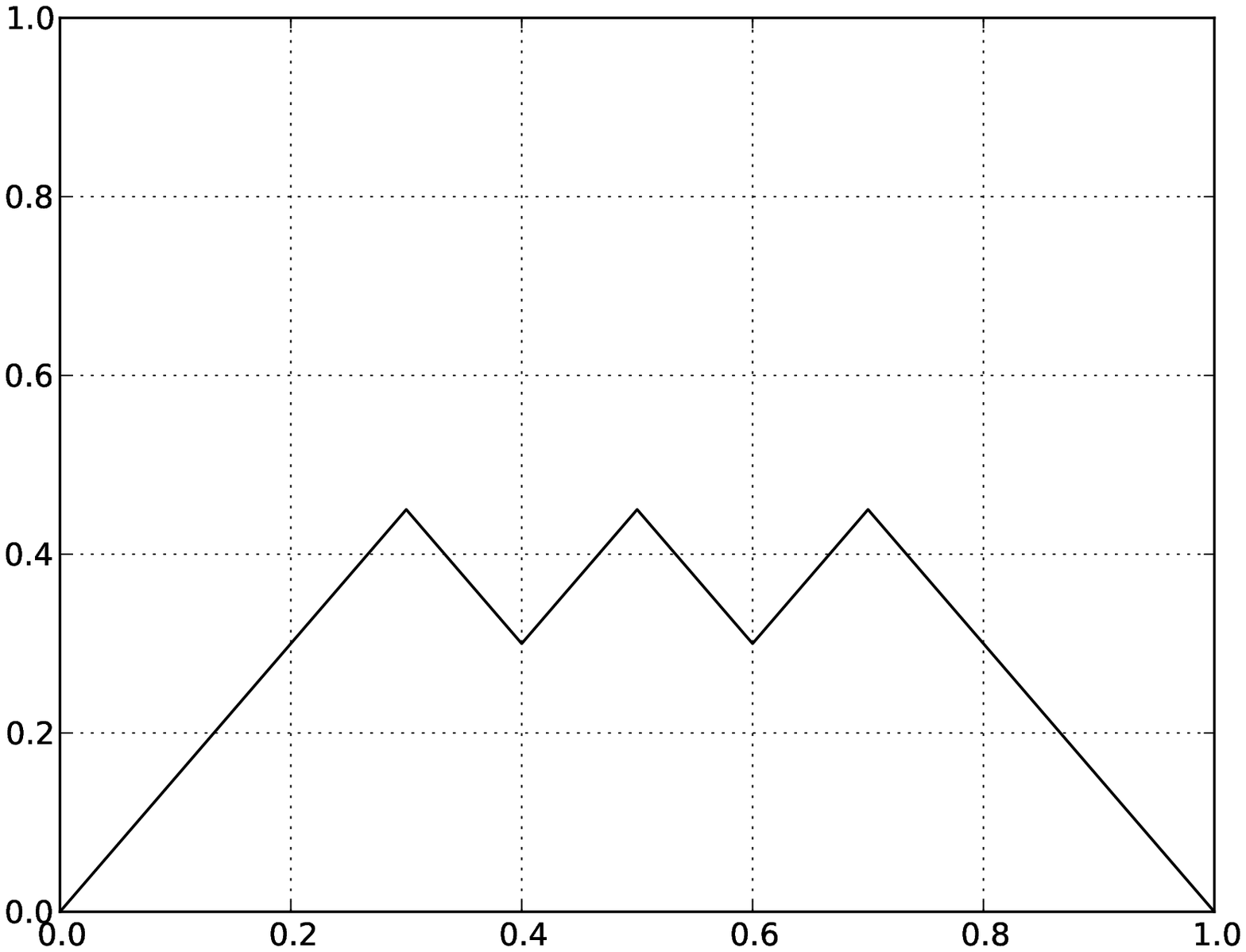}
\caption{A piecewise linear 4-modal map (left) and 5-modal map (right) with
constant slope $s=\pm1.5$.}
\label{fig7}
\end{figure}

Table 4 summarizes the computational performance of both algorithms with the 
$4$-modal map. As happened with the 2-, and 3-modal maps (Tables 2 and 3),
the new algorithm needs more computation loops but less execution time than
the old one for all $\varepsilon=10^{-d}$, $4\leq d\leq7$. 
\begin{table}[h]
\centering
\begin{tabular}{|c|c|c|c|c|c|c|}
\hline
& $h_{old}$ & $n_{old}$ & $t_{old}$ & $h_{new}$ & $n_{new}$ & $t_{new}$ \\ 
\hline
$\varepsilon=10^{-4}$ & 0.421218 & 160 & 0.697776 & 0.422215 & 169 & 0.576800
\\ \hline
$\varepsilon=10^{-5}$ & 0.410476 & 503 & 6.583444 & 0.410776 & 533 & 5.747668
\\ \hline
$\varepsilon=10^{-6}$ & 0.407051 & 1589 & 65.236068 & 0.407147 & 1683 & 
57.951979 \\ \hline
$\varepsilon=10^{-7}$ & 0.405967 & 5021 & 678.706894 & 0.405997 & 5321 & 
616.59469 \\ \hline
\end{tabular}%
\caption{Comparison of performances when computing $h(f)$ in nats with the $%
4 $-modal map of Fig. 7 (left).}
\label{table4}
\end{table}

Likewise, Table 5 summarizes the computational performance of both
algorithms with the $5$-modal map. It is worth noting that now both
algorithms need the same number of loops for all halt criteria $\varepsilon$%
, and yet the new algorithm is faster. 
\begin{table}[h]
\centering
\begin{tabular}{|c|c|c|c|c|c|c|}
\hline
& $h_{old}$ & $n_{old}$ & $t_{old}$ & $h_{new}$ & $n_{new}$ & $t_{new}$ \\ 
\hline
$\varepsilon=10^{-4}$ & 0.420542 & 152 & 0.848166 & 0.420542 & 152 & 0.644305
\\ \hline
$\varepsilon=10^{-5}$ & 0.410239 & 480 & 8.231152 & 0.410239 & 480 & 6.501711
\\ \hline
$\varepsilon=10^{-6}$ & 0.406978 & 1515 & 81.429872 & 0.406978 & 1515 & 
65.307619 \\ \hline
$\varepsilon=10^{-7}$ & 0.405944 & 4788 & 864.376277 & 0.405944 & 4788 & 
695.24749 \\ \hline
\end{tabular}%
\caption{Comparison of performances when computing $h(f)$ in nats with the $%
5 $-modal map of Fig. 7 (right).}
\label{table5}
\end{table}

As in the preceding simulations, we conclude from Table 4 and 5 that both
algorithms determine two correct decimal positions of the topological
entropy of the corresponding map, $h(f)=0.40...$ nats. But this time the
halt criterion $\varepsilon=10^{-6}$ does not suffice; here one has to set $%
\varepsilon=10^{-7}$ to achieve the same precision.

A concluding observation. As anticipated in the remark R2 of Sect. 5 and
illustrated in the Tables 1-5, the values of $h_{new}$ converge from above
with ever more computation loops (or smaller values of the parameter $%
\varepsilon$). This property follows for $h_{old}$ from \cite[Thm. 4.2.4]%
{Alseda}.

\section{Conclusion}

The main contributions of this paper are the following.

(i) In Theorem 1, we proved that the transition rules for min-max symbols (%
\ref{transition}) and (\ref{transition2}), which were derived in \cite%
{Amigo2012} for twice differentiable multimodal maps, actually hold true for
just continuous ones.

(ii) As a result of Theorem 1, we conclude that the validity of formula (\ref%
{sinu2}), which was proved in \cite[Theorem 5.3]{Amigo2012} for twice
differentiable multimodal maps, can be extended to continuous maps. For
subsequent applications, only the particularization of (\ref{sinu2}) to
boundary-anchored maps (Theorem 2) is needed.

(iii) The results reviewed and proved in Sects. 2 and 3, leads to the closed
formula (\ref{efficient}) for the topological entropy of multimodal maps.
Previously we proved in Theorem 3 that, although $\ell_{n}$ clearly depends
on the boundary conditions, the limit $h(f)=\lim_{n\rightarrow\infty}\frac{1%
}{n}\log\ell_{n}$ does not.

(iv) The numerical algorithm proposed in Sect. 5 for the computation of $%
h(f) $ amounts to a recursive scheme to approximate the limit in the closed
formula (\ref{efficient}).

This algorithm is a simplification and, at the same time, a generalization
of the recursion scheme proposed in \cite{Amigo2012} for $h(f)$. Indeed, it
is a simplification because Eqn. (\ref{sinu2}) was used in \cite{Amigo2012}
to compute the lap number $\ell_{\nu}$, while the abridged expression (\ref%
{sinu}) is used here. In other words, the new algorithm does not track the
orbits of the endpoints. And it is also a generalization because we proved
in Theorem 2 that (\ref{sinu}) (and (\ref{sinu2}) for that matter) holds not
only for twice differentiable maps (as assumed in \cite[Theorem 5.3]%
{Amigo2012})) but also for just continuous ones. By the way, this point was
numerically checked in Sect. 6.4.

The performances of both algorithms, old and new, were compared in Sect. 6.1
to 6.4 using smooth and non-smooth $l$-modal maps with $1\leq l\leq5$. In
view of the results summarized in Tables 1 to 5, the old algorithm performs
better in the unimodal case, while the opposite occurs in the other
multimodal cases.

\section{Acknowledgements}

We thank our referees for their constructive criticism.
We are also grateful to Jos\'{e} S. C\'{a}novas and
Mar\'{\i}a Mu\~{n}oz Guillermo (Universidad Polit\'{e}cnica de Cartagena,
Spain) for clarifying discussions, and to V\'{\i}ctor Jim\'{e}nez (Universidad
de Murcia, Spain) for the elegant proof in the Appendix. This work was
financially supported by the Spanish \textit{Ministerio de Econom\'{\i}a y
Competitividad}, grant MTM2012-31698.


\bibliographystyle{unsrt}
\bibliography{SIMPLIFIEDALGOarXiv2022}


%
%
%
%

\section{APPENDIX}

Let $g:X\rightarrow X$ be a continuous map of a compact Hausdorff space $X$
into itself. A point $x\in X$ is nonwandering with respect to the map $g$ if
for any neighborhood $U$ of $x$ there is an $n\geq 1$ (possibly depending on 
$x$) such that $f^{n}(U)\cap U\neq \emptyset $. Fixed and periodic points
are examples of nonwandering points. The closed set of all nonwandering
points of $g$ is called its \textit{nonwandering set} and denoted by $\Omega
(g)$. According to \cite[Lemma 4.1.5]{Alseda}, 
\begin{equation}
h(g)=h(\left. g\right\vert _{\Omega (g)}).  \label{fact1}
\end{equation}%
Furthermore, if 
\begin{equation*}
X=\bigcup\limits_{i=1}^{k}Y_{i}
\end{equation*}%
and all $Y_{i}$ are closed and $g$-invariant (i.e., $g(Y_{i})\subset Y_{i})$%
, then \cite[Lemma 4.1.10]{Alseda},%
\begin{equation}
h(g)=\max_{1\leq i\leq k}h(\left. g\right\vert _{Y_{i}}).  \label{fact2}
\end{equation}

To prove Theorem 3, suppose that $f$ is an $l$-modal selfmap of the compact
interval $I$ with positive shape (the proof for maps with negative shape is
analogous).

Set $I=[a,b]$, and $J=[a^{\prime},b^{\prime}]$ with $a^{\prime}\leq a<b\leq
b^{\prime}$. If $f(a)=a$, choose $a^{\prime}=a$; if $f(b)=a$ ($l$ odd) or $%
f(b)=b$ ($l$ even), choose $b^{\prime}=b$. For definiteness, we suppose the
most general situation, namely, $a^{\prime}<a$ \ and $b<b^{\prime}$. Let $%
F:J\rightarrow J$ be such that (i) $F$ is strictly increasing on $[a^{\prime
},a]$, (ii) $\left. F\right\vert _{[a,b]}=f$, and (iii) $F$ is strictly
decreasing ($l$ odd) or strictly increasing ($l$ even) on $[b,b^{\prime}]$.
In particular, $F$ may be taken piecewise linear on $[a^{\prime},a]\cup%
\lbrack b,b^{\prime}]$. Thus, $F\in\mathcal{M}_{l}(J)$ has the same critical
points and values as $f$, has the same shape and is boundary-anchored. Note
that the shape enters in how $f$ is extended to $F$.

Moreover, it is easy to check that $\Omega (F)=\Omega (f)\cup C$, where $C$
is a closed and $F$-invariant set that only contains fixed points. Thus, $%
h(\left. F\right\vert _{C})=0$ and, according to (\ref{fact1}) and (\ref%
{fact2}), 
\begin{equation*}
h(F)=h(\left. F\right\vert _{\Omega (F)})=\max \{h(\left. F\right\vert
_{\Omega (f)}),h(\left. F\right\vert _{C})\}=h(\left. F\right\vert _{\Omega
(f)})=h(\left. f\right\vert _{\Omega (f)})=h(f).\;\square
\end{equation*}

\end{document}